\documentstyle[11pt]{article}
\newlength{\standardunitlength}
\newtheorem{cor}{Corollary} 
\newtheorem{theorem}{Theorem} 
\newenvironment{proof}{\noindent {\sc Proof:}}{$\Box$ \vspace{2 ex}}

\begin{document}

\begin{center} {\bf A Probabilistic Proof of the Rogers-Ramanujan
Identities} \end{center}

\begin{center}
By Jason Fulman
\end{center}

\begin{center}
Stanford University
\end{center}

\begin{center}
Department of Mathematics
\end{center}

\begin{center}
Building 380, MC 2125
\end{center}

\begin{center}
Stanford, CA 94305, USA
\end{center}

\begin{center}
fulman@math.stanford.edu
\end{center}

\begin{center}
http://math.stanford.edu/$\sim$ fulman
\end{center}

\begin{center}
Math Reviews Subject Classification: 20P05
\end{center}

\newpage

\begin{abstract} The asymptotic probability theory of conjugacy
classes of the finite general groups leads to a probability measure on
the set of all partitions of natural numbers. A simple method of
understanding these measures in terms of Markov chains is given,
leading to an elementary probabilistic proof of the Rogers-Ramanujan
identities. This is compared with work on the uniform measure. The
main case of Bailey's lemma is interpreted as finding eigenvectors of
the transition matrix of a Markov chain. It is shown that the
viewpoint of Markov chains extends to quivers. \end{abstract}

\begin{center}
\end{center}

\section{Introduction}

	The Rogers-Ramanujan identities \cite{Ro}, \cite{RoRa},
\cite{Sc}

\[ 1+\sum_{n=1}^{\infty} \frac{q^{n^2}}{(1-q)(1-q^2) \cdots (1-q^n)} =
\prod_{n=1}^{\infty} \frac{1}{(1-q^{5n-1})(1-q^{5n-4})} \]

\[ 1+\sum_{n=1}^{\infty} \frac{q^{n(n+1)}}{(1-q)(1-q^2) \cdots
(1-q^n)} = \prod_{n=1}^{\infty} \frac{1}{(1-q^{5n-2})(1-q^{5n-3})} \]
are among the most interesting partition identities in number theory
and combinatorics. Combinatorial aspects of these identities are
discussed by Andrews \cite{A},\cite{A2}. Lepowsky and Wilson
\cite{LP1}, \cite{LP2}, \cite{LP3} connect the Rogers-Ramanujan
identities with affine Lie algebras and conformal field theory. Feigen
and Frenkel \cite{FF} interpret them as a character formula for the
Virasoro algebra. The author \cite{F2} applies the product form to
computational group theory, obtaining a simple formula for the $n
\rightarrow \infty$ limit of the probability that an element of
$GL(n,q)$ is semisimple. Stembridge \cite{St} proves them by adapting
a method of Macdonald for calculating partial fraction expansions of
symmetric formal power series. Garsia and Milne \cite{GM} offer a lengthy
bijective proof, giving birth to the ``involution principle'' in
combinatorics. Garrett, Ismail, and Stanton \cite{GIT} usefully recast
Roger's original proof using orthogonal polynomials. Relations with
statistical mechanics appear in \cite{Ba1}, \cite{ABF}, and
\cite{W}. The paper \cite {BeMc} gives an account of their appearance
in physics.

	One of the main results of this note is a first probabilistic
proof of the Rogers-Ramanujan identities. In fact the $i=1$ and $i=k$
cases of Theorem \ref{Gordon} will be proved. Theorem \ref{Gordon} is due
to Andrews \cite{A4}. Bressoud \cite{Br3} then connected it with an
earlier combinatorial result of Gordon \cite{Go}.

\begin{theorem} \label{Gordon} (\cite{A4}) Let $x_n$ denote
$(1-x)\cdots(1-x^n)$. For $1\leq i \leq k, k \geq 2$,

\[ \sum_{n_1,\cdots,n_{k-1} \geq 0} \frac{x^{N_1^2 + \cdots +
N_{k-1}^2+ N_i + \cdots+ N_{k-1}}}{(x)_{n_1}\cdots(x)_{n_{k-1}}} =
\prod_{r=1 \atop r \neq 0, \pm i (mod \ 2k+1)}^{\infty}
\frac{1}{1-x^r} \] where $N_j = n_j + \cdots n_{k-1}$. \end{theorem}
There are other elementary proofs of the Rogers-Ramanujan identities
available (e.g. \cite{Br1}), but the one offered here is possibly the
most natural and gives insight into Bailey's Lemma.

	The basic object of study in this paper is a certain
one-parameter family of probability measures $M_u$ on the set of all
partitions of all natural numbers, studied in the prior article
\cite{F1}. Section \ref{measure} recalls these measures and gives
their group theoretic motivation. (For now we remark that they arise
in the study of conjugacy classes of finite classical groups. As
conjugacy classes of compact Lie groups are essentially eigenvalues up
to the action of the Weyl group, their probabilistic study can be
regarded as philosophically analogous to work of Dyson \cite{Dy}, who
described the eigenvalues of random matrices of compact Lie groups in
terms of Brownian motion).

	Section \ref{markov} shows how to construct these measures
using a non-upward-moving Markov chain on the integers. The Markov
chain is diagonalizable with eigenvalues $1,\frac{u}{q},
\frac{u^2}{q^4}, \cdots$, and a basis of eigenvectors is
given. Analogous computations are done for a measure related to the
uniform measure on partitions, which by work of Fristedt \cite{Fr} has
a Markov chain approach. It would be interesting to make a connection
with the articles \cite{VS}, in which a fascinating continuous space
Markov chain arises in the asymptotic probability theory of the
symmetric group.

	With these preliminaries in place, Section \ref{RogersRam}
gives a proof of the Rogers-Ramanujan identities. The idea of the
proof is simple. We compute in two ways the $L \rightarrow \infty$
probability that the Markov chain started at $L$ is absorbed at the
point $0$ after $k$ steps. (Since the Markov chain is absorbed at 0
with probability 1 and the measure $M_u$ corresponds to the
$L \rightarrow \infty$ limit, the time to absorption really is the most
natural quantity one could consider). The sum side of the
Rogers-Ramanujan identities follows from the definition of the
probability measures. For the product side, the fact that the
transition matrix is explicitly diagonalizable gives a sum
expression. One then applies the Jacobi triple product identity (which
as explained on page 21 of \cite{A} follows easily from the
$q$-binomial theorem) and the proof is complete.

	Section \ref{RogersRam} continues by discussing Bailey's
Lemma, which is the only non-trivial step in many of the simplest
proofs of the Rogers-Ramanujan identities. The most useful case of
Bailey's Lemma follows immediately once a basis of eigenvectors of the
transition matrix has been found. The only non-trivial step in our
proof of the Rogers-Ramanujan identities is finding a basis of
eigenvectors; however with Mathematica this was easy. By contrast it
is unclear how one would guess at Bailey's Lemma.

	Section \ref{quivers} reviews the theory of quivers and shows
that the Markov chain method extends to quivers. Although we have not
invested serious effort into finding analogs of Bailey's Lemma for
quivers other than the one point quiver (which corresponds to
conjugacy classes of the finite general linear groups), it is not hard
to see that the resulting Bailey Lemmas differ from those of \cite{ML}
and \cite{ASW}.

	It is tempting to speculate that there is a relationship
between conjugacy classes of the finite general linear groups and
modular forms. Aside from this paper, there are two good reasons to
suspect this. One reason is that the conjugacy classes are related to
Hall-Littlewood polynomials, which in turn are related to vertex
operators \cite{J}. Other evidence is work of Bloch and Okounkov
\cite{BO}, who relate a version of the uniform measure on partitions
to quasi-modular forms.

\section{Measures on Partitions and Group Theory} \label{measure}

	We begin by reviewing some standard notation about partitions,
as on pages 2-5 of Macdonald $\cite{Mac}$. Let $\lambda$ be a
partition of some non-negative integer $|\lambda|$ into parts
$\lambda_1 \geq \lambda_2 \geq \cdots$. Let $m_i(\lambda)$ be the
number of parts of $\lambda$ of size $i$, and let $\lambda'$ be the
partition dual to $\lambda$ in the sense that $\lambda_i' =
m_i(\lambda) + m_{i+1}(\lambda) + \cdots$. Let $n(\lambda)$ be the
quantity $\sum_{i \geq 1} (i-1) \lambda_i$. It is also useful to
define the diagram associated to $\lambda$ as the set of points $(i,j)
\in Z^2$ such that $1 \leq j \leq \lambda_i$. We use the convention
that the row index $i$ increases as one goes downward and the column
index $j$ increases as one goes across. So the diagram of the
partition $(5441)$ is:

\[ \begin{array}{c c c c c}
		. & . & . & . & .  \\
		. & . & . & . &    \\
		. & . & . & . &    \\
		. & & & &  
	  \end{array}  \]

	The rest of this section follows the paper \cite{F1}. Let $q$
be the size of a finite field. To begin we recall a way of defining a
one parameter family of probability measures $M_u(\lambda)$ on
the set of all partitions of all natural numbers. If one simply wants
a formula, then all of the following definitions are equivalent. In
the third expression, $P_{\lambda}$ denotes a Hall Littlewood
polynomial as in \cite{Mac}.

\begin{eqnarray*}
M_u(\lambda) & = & \left[\prod_{r=1}^{\infty} (1-\frac{u}{q^r})\right]
\frac{u^{|\lambda|}} {q^{\sum_i (\lambda'_i)^2} \prod_{i} (\frac{1}{q})_{m_i(\lambda)}}\\
& = & \left[\prod_{r=1}^{\infty} (1-\frac{u}{q^r})\right]
\frac{u^{|\lambda|}}{q^{2 [\sum_{h<i} h m_h(\lambda) m_i(\lambda) +
\frac{1}{2} \sum_i (i-1) m_i(\lambda)^2]}\prod_i |GL(m_i(\lambda),q)|}\\
& = & \left[\prod_{r=1}^{\infty} (1-\frac{u}{q^r})\right]
\frac{u^{|\lambda|} P_{\lambda}(\frac{1}{q},\frac{1}{q^2},\cdots;0,\frac{1}{q})}{q^{
n(\lambda)}}
\end{eqnarray*}

	The only fact from this section which is required for the
proof of the Rogers-Ramanujan identities is the fact that $M_u$
defines a probability measure for $0<u<1$. This can be seen using an
identity from either \cite{Mac} or else following Stong \cite{Sto},
who uses the fact that there are $q^{n^2-n}$ unipotent elements in
$GL(n,q)$. The first proof of Theorem \ref{markovgl} will use the fact
that $M_u$ is a probability measure without further comment. As there
has been interest in simplifying the proofs of the Rogers-Ramanujan
identities as much as possible, a second completely elementary proof
of Theorem \ref{markovgl} will be given. From this second proof it
will follow that $M_u(\lambda)$ is a probability measure.

	Although not logically necessary for this paper, we mention
that for $0<u<1$ and $q$ a prime power, the measures $M_u$ have a
group theoretic description, referring the reader to the survey
\cite{F4} for further discussion. Recall that the conjugacy
classes of $GL(n,q)$ are parameterized by rational canonical
form. Each such matrix corresponds to the following combinatorial
data. To every monic non-constant irreducible polynomial $\phi$ over
$F_q$, associate a partition (perhaps the trivial partition)
$\lambda_{\phi}$ of some non-negative integer $|\lambda_{\phi}|$. The
only restrictions necessary for this data to represent a conjugacy
class are that $|\lambda_z| = 0$ and $\sum_{\phi} |\lambda_{\phi}|
deg(\phi) = n.$ To be explicit, a representative of the conjugacy
class corresponding to the data $\lambda_{\phi}$ may be given as
follows. Define the companion matrix $C(\phi)$ of a polynomial
$\phi(z)=z^{deg(\phi)} + \alpha_{deg(\phi)-1} z^{deg(\phi)-1} + \cdots
+ \alpha_1 z + \alpha_0$ to be:

\[ \left( \begin{array}{c c c c c}
                0 & 1 & 0 & \cdots & 0 \\
                0 & 0 & 1 & \cdots & 0 \\
                \cdots & \cdots & \cdots & \cdots & \cdots \\
                0 & 0 & 0 & \cdots & 1 \\
                -\alpha_0 & -\alpha_1 & \cdots & \cdots & -\alpha_{deg(\phi)-1}
          \end{array} \right) \] Let $\phi_1,\cdots,\phi_k$ be the
polynomials such that $|\lambda_{\phi_i}|>0$. Denote the parts of
$\lambda_{\phi_i}$ by $\lambda_{\phi_i,1} \geq \lambda_{\phi_i,2} \geq
\cdots $. Then a matrix corresponding to the above conjugacy class
data is:

\[ \left( \begin{array}{c c c c}
 C(\phi_i^{\lambda_{\phi_i,1}}) & 0 & 0  \\
                0 & C(\phi_i^{\lambda_{\phi_i,2}}) & 0 \\
                0 & 0 & \cdots 
          \end{array} \right) \]
        
        Now consider the following procedure for putting a measure on
the set of all partitions of all natural numbers. Fix $u$ such that
$0<u<1$. Pick a non-negative integer such that the chance of choosing
$n$ is equal to $(1-u)u^n$. Then pick $\alpha$ uniformly in $GL(n,q)$
and take $\lambda$ to be the paritition corresponding to the
polynomial $z-1$ in the rational canonical form of $\alpha$. If $n=0$
take $\lambda$ to be the trivial partition. The random partition so
defined obeys $M_u$ measure.  (The polynomial $z-1$ is considered
without loss of generality. Partitions corresponding to other
irreducible polynomials are probabilistically independent, and one
just replaces $q$ by raising it to the degree of the polynomial). In
the limit $u \rightarrow 1$, one is simply studying random elements in
a fixed $GL(n,q)$ with $n \rightarrow \infty$.  The substitutions $u
\rightarrow -u$ and $q \rightarrow -q$ correspond to the finite
unitary groups. The idea of auxilliary randomization of the dimension
$n$ is analogous to the idea of canonical ensembles in statistical
mechanics.

\section{Markov chains} \label{markov}

\subsection{Group theoretical measures}

	The first result of this paper describes the measure $M_u$ in
terms of Markov chains. Two proofs will be given. The first proof is
given in the interest of clarity and assumes that $M_u$ is a
probability measure. The second proof is more elementary.

	It is convenient to set $\lambda_0'$ (the height of an
imaginary zeroth column) equal to $\infty$. For the entirety of this
subsection, $(x)_n$ will denote $(1-x)(1-x/q) \cdots
(1-x/q^{n-1})$. Thus $(x)_0=1$ and $(x)_n=0$ for $n<0$. For
convenience of notation, let $P(a)$ be the $M_u$ probability that
$\lambda_1=a$. $Prob(E)$ will denote the probability of an event $E$
under the measure $M_u$.

\begin{theorem} \label{markovgl} Starting with $\lambda_0'=\infty$,
define in succession $\lambda_1',\lambda_2',\cdots$ according to the
rule that if $\lambda_i'=a$, then $\lambda_{i+1}'=b$ with probability

\[ K(a,b) = \frac{u^b (\frac{1}{q})_a (\frac{u}{q})_a}{q^{b^2}
(\frac{1}{q})_{a-b} (\frac{1}{q})_b (\frac{u}{q})_b}.\] Then the
resulting partition is distributed according to $M_u$. \end{theorem}

\begin{proof} The $M_u$ probability of choosing a partition with
$\lambda_i'=r_i$ for all $i$ is

\[ Prob.(\lambda_0'=\infty) \frac{Prob.(\lambda_0'=\infty,
\lambda_1'=r_1)} {Prob.(\lambda_0'=\infty)} \prod_{i=1}^{\infty}
\frac{Prob.(\lambda_0'=\infty,
\lambda_1'=r_1,\cdots,\lambda_{i+1}'=r_{i+1})}
{Prob.(\lambda_0'=\infty,\lambda_1'=
r_1,\cdots,\lambda_{i}'=r_{i})}.\] Thus it is enough to prove that

\[  \frac{Prob.(\lambda_0'=\infty,
\lambda_1'=r_1,\cdots,\lambda_{i-1}'=r_{i-1}, \lambda_i'=a,
\lambda_{i+1}'=b)} {Prob.(\lambda_0'=\infty,\lambda_1'=
r_1,\cdots,\lambda_{i-1}'=r_{i-1}, \lambda_{i}'=a)} = \frac{u^b
(\frac{1}{q})_a (\frac{u}{q})_a}{q^{b^2} (\frac{1}{q})_{a-b}
(\frac{1}{q})_b (\frac{u}{q})_b},\] for all
$i,a,b,r_1,\cdots,r_{i-1} \geq 0$.

	For the case $i=0$, the equation \[ P(a) = \frac{u^a
(\frac{u}{q})_{\infty}}{q^{a^2} (\frac{1}{q})_a (\frac{u}{q})_a} \] is
given a probabilistic proof in \cite{F1}. For an elementary proof of
this identity, see the second proof of this theorem. For $i>0$ one
calculates that

\[ \sum_{\lambda: \lambda_1'=r_1,\cdots,\lambda_{i-1}'=r_{i-1} \atop
\lambda_i'=a} M_u(\lambda) = \frac{u^{r_1+\cdots+r_{i-1}}}
{q^{r_1^2+\cdots+r_{i-1}^2} (\frac{1}{q})_{r_1-r_2} \cdots
(\frac{1}{q})_{r_{i-2}-r_{i-1}} (\frac{1}{q})_{r_{i-1}-a}} P(a).\]
Similarly, observe that \[ \sum_{\lambda:
\lambda_1'=r_1,\cdots,\lambda_{i-1}'=r_{i-1} \atop
\lambda_i'=a,\lambda_{i+1}'=b} M_u(\lambda) =
\frac{u^{r_1+\cdots+r_{i-1}+a}} {q^{r_1^2+\cdots+r_{i-1}^2+a^2}
(\frac{1}{q})_{r_1-r_2} \cdots (\frac{1}{q})_{r_{i-2}-r_{i-1}}
(\frac{1}{q})_{r_{i-1}-a} (\frac{1}{q})_{a-b}} P(b).\] Thus the ratio
of these two expressions is \[ \frac{u^b (\frac{1}{q})_a
(\frac{u}{q})_a}{q^{b^2} (\frac{1}{q})_{a-b} (\frac{1}{q})_b
(\frac{u}{q})_b}, \] as desired. Note that the transition
probabilities automatically sum to 1 because

\[ \sum_{b \leq a} \frac{\sum_{\lambda:
\lambda_1'=r_1,\cdots,\lambda_{i-1}'=r_{i-1} \atop
\lambda_i'=a,\lambda_{i+1}'=b} M_u(\lambda)}{\sum_{\lambda:
\lambda_1'=r_1,\cdots,\lambda_{i-1}'=r_{i-1} \atop \lambda_i'=a}
M_u(\lambda)} =1.\]

\end{proof}

\begin{proof} (Second Proof) This proof needs only that $M_u$ is a
measure; it will emerge that $M_u$ is a probability measure. For this
proof $P(a)$ denotes the $M_u$ mass that $\lambda_1'=a$.

	One calculates that

\[ \sum_{\lambda: \lambda_1'=r_1,\cdots,\lambda_{i-1}'=r_{i-1} \atop
\lambda_i'=a} M_u(\lambda) = \frac{u^{r_1+\cdots+r_{i-1}}}
{q^{r_1^2+\cdots+r_{i-1}^2} (\frac{1}{q})_{r_1-r_2} \cdots
(\frac{1}{q})_{r_{i-2}-r_{i-1}} (\frac{1}{q})_{r_{i-1}-a}} P(a).\]
Similarly, observe that \[ \sum_{\lambda:
\lambda_1'=r_1,\cdots,\lambda_{i-1}'=r_{i-1} \atop
\lambda_i'=a,\lambda_{i+1}'=b} M_u(\lambda) =
\frac{u^{r_1+\cdots+r_{i-1}+a}} {q^{r_1^2+\cdots+r_{i-1}^2+a^2}
(\frac{1}{q})_{r_1-r_2} \cdots (\frac{1}{q})_{r_{i-2}-r_{i-1}}
(\frac{1}{q})_{r_{i-1}-a} (\frac{1}{q})_{a-b}} P(b).\] Thus the ratio
of these two expressions is \[ \frac{P(b)u^a}{P(a) q^{a^2}
(\frac{1}{q})_{a-b}}.\] Since $M_u$ is a measure, it follows that \[
\sum_{b \leq a} \frac{P(b)u^a}{P(a) q^{a^2} (\frac{1}{q})_{a-b}}=1.\]
From this recursion and the fact that $P(0)=(\frac{u}{q})_{\infty}$,
one solves for $P(a)$ inductively, finding that \[ P(a) = \frac{u^a
(\frac{u}{q})_{\infty}} {q^{a^2} (\frac{1}{q})_a (\frac{u}{q})_a} .\]
Identity 2.2.8 on page 20 of \cite{A} now implies that $\sum_a
P(a)=1$, so that $M_u$ is a probability measure.  \end{proof}
	
	The algorithm of Theorem \ref{markovgl} is nice in that it can
be implemented on a computer. An analogous Markov chain approach to
the conjugacy classes of the finite symplectic and orthogonal groups
has recently been worked out in \cite{F3}.

	Theorem \ref{diagonalize} explicitly diagonalizes the
transition matrix $K$, which is fundamental for understanding the
Markov chain. Note that if the current distribution of the Markov
chain is given by a row vector, the distribution at the next step is
obtained by multiplying the row vector on the right by $K$.

\begin{theorem} \label{diagonalize} \begin{enumerate}

\item Let $C$ be the diagonal matrix with $(i,i)$ entry
$(\frac{1}{q})_i (\frac{u}{q})_i$. Let $M$ be the matrix $\left(
\frac{u^j}{q^{j^2} (\frac{1}{q})_{i-j}} \right)$. Then $K=CMC^{-1}$,
which reduces the problem of diagonalizing $K$ to that of
diagonalizing $M$.

\item Let $A$ be the matrix $\left( \frac{1}{(\frac{1}{q})_{i-j}
(\frac{u}{q})_{i+j}} \right)$. Then the columns of $A$ are
eigenvectors of $M$ for right multiplication, the $j$th column having
eigenvalue $\frac{u^j}{q^{j^2}}$.

\item The inverse matrix $A^{-1}$ is $\left( \frac{(1-u/q^{2i})
(-1)^{i-j} (\frac{u}{q})_{i+j-1}}{q^{i-j \choose 2}
(\frac{1}{q})_{i-j}} \right)$.

\end{enumerate}
\end{theorem}

\begin{proof} The first part is obvious. The second part is a special
case of Lemma 1 of \cite{Br1}. The third part is a lemma of
\cite{An24}. \end{proof}

	The point of the proof of Theorem \ref{diagonalize} is that
once one knows (either from Mathematica or by implementing algorithms
from linear algebra) what the eigenvectors are, it is a simple matter
to verify the computation.

	The following corollary will be useful for the proof of the
Rogers-Ramanujan identities.

\begin{cor} \label{bigcor} Let $E$ be the diagonal matrix with $(i,i)$
entry $\frac{u^i}{q^{i^2}}$. Then $K^r=C A E^r A^{-1} C^{-1}$. More
explicitly,

\[ K^r(L,j) = \frac{(\frac{1}{q})_L (\frac{u}{q})_L}{(\frac{1}{q})_j
(\frac{u}{q})_j} \sum_{n=0}^{\infty} \frac{u^{rn} (1-u/q^{2n})
(-1)^{n-j} (\frac{u}{q})_{n+j-1}}{q^{rn^2} (\frac{1}{q})_{L-n}
(\frac{u}{q})_{L+n} q^{n-j \choose 2} (\frac{1}{q})_{n-j}}.\] \end{cor}

\begin{proof} This is immediate from Theorem
\ref{diagonalize}. \end{proof}

\subsection{Mixture of uniform measures}

	For this subsection $q<1$. The measure assigns probability
$q^{|\lambda|} \prod_{i=1}^{\infty} (1-q^i)$ to the partition
$\lambda$. Conditioning this measure to live on partitions of a given
size gives the uniform measure, an observation exploited by Fristedt
\cite{Fr}. As is clear from \cite{O}, this measure is very natural
from the viewpoint of representation theory.

	Fristedt (loc. cit.) proved that this measure has a Markov
chain description. His chain affects row lengths rather than column
lengths (though the algorithm would work on columns too as the measure
is invariant under transposing diagrams). Nevertheless, we adhere to
his notation. We use the notation that $(x)_n=(1-x)\cdots(1-x^n)$.

\begin{theorem} \cite{Fr}  Starting with $\lambda_0=\infty$,
define in succession $\lambda_1,\lambda_2,\cdots$ according to the
rule that if $\lambda_i=a$, then $\lambda_{i+1}=b$ with probability

\[ K(a,b) = \frac{q^b (q)_a}{(q)_b}.\] Then the resulting partition is
distributed according to the measure of this subsection. \end{theorem}

	Theorem \ref{diag2} diagonalizes this Markov chain, giving a
basis of eigenvectors. The proof is analogous to that of Theorem
\ref{diagonalize}, the second part being proved by induction.

\begin{theorem} \label{diag2}

\begin{enumerate}

\item Let $C$ be the diagonal matrix with $(i,i)$ entry equal to
$\frac{(q)_i}{q^i}$. Let $M$ be the matrix with $(i,j)$ entry $q^i$ if $i
\geq j$ and $0$ otherwise. Then $K=CMC^{-1}$.

\item Let $D$ be the diagonal matrix with $(i,i)$ entry equal to
$q^i$. Let $A$ be the matrix $\left( \frac{(-1)^{i-j}}{q^{i-j \choose
2} (\frac{1}{q})_{i-j}} \right)$, so that its inverse is
$\left(\frac{1}{(\frac{1}{q})_{i-j}} \right)$ by part b of Theorem
\ref{diagonalize}. Then the eigenvectors of $M$ are the columns of
$A$, the $j$th column having eigenvalue $q^j$.

\end{enumerate}
\end{theorem}

	As a corollary, one obtains a simple expression for the chance
that under the measure of this subsection, the $r$th row has size $j$.

\begin{cor} \[ K^r(L,j) = \frac{q^j q^{L(r-1)} (q)_L
(\frac{1}{q})_{L-j+r-1}}{(q)_j (\frac{1}{q})_{L-j}
(\frac{1}{q})_{r-1}}.\] Letting $L \mapsto \infty$, the chance that
the $r$th row has size $j$ becomes \[ \frac{(q)_{\infty} q^{rj}}{(q)_j
(q)_{r-1}}.\] \end{cor}

\begin{proof} To obtain the first expression, one multiplies out
$K^r=CAD^rA^{-1}C^{-1}$ and uses the $q$-binomial theorem 

\[ \sum_{m=0}^{\infty} y^m q^{(m^2+m)/2} \frac{(q)_n}{(q)_m (q)_{m-n}}
= (1+yq)(1+yq^2) \cdots (1+yq^n).\] \end{proof}

	The second part of the corollary can be proved directly
without recourse to Markov chain theory; one simply attaches to an
$r*j$ square two partitions: one with at most $r-1$ rows and another
with at most $j$ columns.

\section{Rogers-Ramanujan Identities and Bailey's Lemma} \label{RogersRam}

	The first result of this section proves the following identity
of Andrews, which contains the Rogers-Ramanujan identities. In this
section $(x)_n$ denotes $(1-x)(1-x/q) \cdots (1-x/q^{n-1})$.

\begin{theorem} \cite{A4} For $k \geq 2$,

\[ \sum_{n_1,\cdots,n_{k-1} \geq 0} \frac{1}{q^{N_1^2 + \cdots +
N_{k-1}^2}(1/q)_{n_1}\cdots(1/q)_{n_{k-1}}} =
\prod_{r=1 \atop r \neq 0, \pm k (mod \ 2k+1)}^{\infty}
\frac{1}{1-(1/q)^r} \]

\[ \sum_{n_1,\cdots,n_{k-1} \geq 0} \frac{1}{q^{N_1^2 + \cdots +
N_{k-1}^2+ N_1 + \cdots+ N_{k-1}}(1/q)_{n_1}\cdots(1/q)_{n_{k-1}}} =
\prod_{r=1 \atop r \neq 0, \pm 1 (mod \ 2k+1)}^{\infty}
\frac{1}{1-(1/q)^r} \]  where $N_j = n_j + \cdots
n_{k-1}$. \end{theorem}

\begin{proof} For the first identity, we compute in two ways \[
\sum_{\lambda: \lambda' <k} M_1(\lambda).\] One obtains the sum side by using the first definition of $M_1$ in Section
\ref{measure}. For the product side, let $u=1$, $j=0$, $r=k$, and $L
\rightarrow \infty$ in Corollary \ref{bigcor}. The rest is now a
standard argument.

\begin{eqnarray*}
1+\sum_{n=1}^{\infty} \frac{(1+1/q^n)(-1)^n}{q^{rn^2} q^{n \choose 2}} & = & 1+\sum_{n=1}^{\infty} (-1)^n (\frac{1}{q})^{(k+1/2)n^2-n/2} + (\frac{1}{q})^{(k+1/2)n^2+n/2}\\
& = & \sum_{n=-\infty}^{\infty} (-1)^n (\frac{1}{q})^{(k+1/2)n^2} (\frac{1}{q})^{n/2}.
\end{eqnarray*} Now invoke Jacobi's triple product identity (e.g. \cite{A})

\[ \sum_{n=-\infty}^{\infty} (-1)^n v^n w^{n^2} =
\prod_{n=1}^{\infty} (1-vw^{2n-1})(1-w^{2n-1}/v)(1-w^{2n}).\] The
proof of the second identity is the same except that one takes
$u=\frac{1}{q}$ instead of $u=1$.
\end{proof}

	Next we discuss the most important case of Bailey's Lemma,
which was alluded to in \cite{Ba} and stated explicitly in
\cite{An24}. A pair of sequences $\{\alpha_L\}$ and $\{\beta_L\}$ are
called a Bailey pair if \[ \beta_L = \sum_{r=0}^L
\frac{\alpha_r}{(1/q)_{L-r} (u/q)_{L+r}}.\] Bailey's Lemma states that
if $\alpha_L'=\frac{u^L}{q^{L^2}} \alpha_L$ and $\beta_L'=\sum_{r=0}^L
\frac{u^r}{q^{r^2}(1/q)_{L-r}} \beta_r$, then $\{\alpha_L'\}$ and
$\{\beta_L'\}$ are a Bailey pair.

	From our viewpoint, this case of Bailey's Lemma is
clear. Namely let $A,D,M$ be as in Theorem \ref{diagonalize} (recall
that $M=ADA^{-1}$). Viewing $\alpha=\vec{\alpha_L}$ and
$\beta=\vec{\beta_L}$ as column vectors, the notion of a Bailey pair
means that $\beta=A \alpha$. This case of Bailey's Lemma follows
because

\[ \beta' = M \beta = ADA^{-1} \beta = AD \alpha = \alpha'.\] It would
be interesting to obtain all of Bailey's Lemma (and Theorem
\ref{Gordon}) by probabilistic arguments, and also to understand the
Bailey lattice probabilistically. The reader may enjoy the survey by
Bressoud \cite{Br2}.

\section{Quivers} \label{quivers}

	This section uses the notion of a quiver, as surveyed in Kac
\cite{K}, to which the reader is referred for more detail. The basic
set-up is as follows. Let $\Gamma$ be a connected graph with $n$
vertices labelled as $1,\cdots,n$ (where we allow loops). Let $N,Z$
denote the natural numbers and integers respectively. Let $f_{ij}$ be
the number of edges between $i,j$. Associated to $\Gamma$ is a natural
bilinear form on $Z^n$ and a root system $\Delta \subset Z^n$. Choose
an arbitrary orientation of $\Gamma$ so that $\Gamma$ is a quiver. For
a given dimension $\alpha \in N^n-\{0\}$, let $A_{\Gamma}(\alpha,q)$
be the number of classes of absolutely indecomposable representations
of $\Gamma$ over the algebraic closure of a field of $q$ elements. It
is proved in \cite{K} that $A_{\Gamma}(\alpha,q)$ is a polynomial in
$q$ with integer coefficients, and that this polynomial is independent
of the orientation of the graph. Kac (loc. cit.) formulated many
conjectures about this polynomial. One such conjecture, which is still
open, is that the constant term in $A_{\Gamma}(\alpha,q)$ is the
multiplicity of $\alpha$ in the root system.

	In recent work, Hua \cite{H1,H2} has given a completely
combinatorial reformulation of this conjecture. To explain, let
$(1/q)_n$ denote $(1-\frac{1}{q}) \cdots (1-\frac{1}{q^n})$, and for
any two partitions $\lambda,\mu$ define $<\lambda,\mu>=\sum_{i \geq 1}
\lambda_i' \mu_i'$. Let $t_j^{\alpha}$ be the coefficient of $q^i$ in
the polynomial $A_{\Gamma}(\alpha,q)$. Let $b_{\lambda}=\prod_{i \geq
1} (\frac{1}{q})_{m_i(\lambda)}$, where $m_i$ is the number of parts
of $\lambda$ of size $i$. Let $\lambda(1),\cdots, \lambda(n)$ be an
$n$-tuple of partitions. Set $U_{\alpha} = U_1^{\alpha_1} \cdots
U_n^{\alpha_n}$. Hua's result, which reduces Kac's conjecture to a
combinatorial assertion, is that

\[ \sum_{\lambda(1),\cdots,\lambda(n)} \frac{\prod_{1 \leq i \leq j
\leq n} q^{f_{ij}<\lambda(i),\lambda(j)>} U_1^{|\lambda(1)|} \cdots U_n^{|\lambda(n)|}}{\prod_{1 \leq i \leq n}
q^{<\lambda(i),\lambda(i)>} b (\lambda(i))} = \prod_{\alpha \in \Delta^+}
\prod_{j=0}^{deg(A_{\Gamma(\alpha,q)})} \prod_{i=1}^{\infty}
(\frac{1}{1-U_{\alpha}q^{j-i}})^{t_j(\alpha)}.\]

	A few points are in order. First, the right hand side of this
equation is different from the statements in \cite{H1,H2}, due to a
minor slip there. Second, observe that the expression converges in the
ring of formal power series in the variables $U_1,\cdots,U_n$. This
leads one to define a ``probability'' measure on $n$-tuples of
partitions $M_{\Gamma,\vec{U}}$ by assigning mass

\[ \prod_{\alpha \in \Delta^+} \prod_{j=0}^{deg(A_{\Gamma(\alpha,q)})}
\prod_{i=1}^{\infty} (1-U_{\alpha}q^{j-i})^{t_j(\alpha)}
\frac{\prod_{1 \leq i \leq j \leq n}
q^{f_{ij}<\lambda(i),\lambda(j)>}}{\prod_{1 \leq i \leq n}
q^{<\lambda(i),\lambda(i)>} b_{\lambda(i)}} \] to the $n$-tuple
$\lambda(1),\cdots,\lambda(n)$. For quivers of finite type this is a
true probability measure for values of $U_1,\cdots,U_n$ sufficiently
small, but in general we abuse notation by using terms from
probability theory when $U_1,\cdots,U_n$ are variables.

	Note that when the graph consists of a single point, this
measure is simply the measure $M_u$ from Section
\ref{measure}. Theorem \ref{stillMarkov} shows that the structure of a
Markov chain is still present. As the idea of the proof is the same as
the second step of Theorem \ref{markovgl}, the algebra is omitted. Let
$P(a_1,\cdots,a_n)$ denote the $M_{\Gamma,\vec{U}}$ probability that
$\lambda(1),\cdots,\lambda(n)$ have $a_1,\cdots,a_n$ parts
respectively.

\begin{theorem} \label{stillMarkov} Let
$\lambda(1)_1',\cdots,\lambda(n)_1'$ be distributed as
$P(a_1,\cdots,a_n)$. Define $(\lambda(1)_2',\cdots,\lambda(n)_2')$
then $(\lambda(1)_3',\cdots,\lambda(n)_3')$, etc. successively
according to the rule that if $(\lambda(1)_i',\cdots,\lambda(n)_i')$
is equal to $(a_1,\cdots,a_n)$, then $(\lambda(1)_{i+1}'
,\cdots,\lambda(n)_{i+1}')$ is equal to $(b_1,\cdots,b_n)$ with
probability

\[ K(\vec{a},\vec{b}) = \prod_{1 \leq i \leq j \leq n} q^{f_{ij} a_i
a_j} \prod_{i=1}^n \frac{U_i^{a_i}}{q^{a_i^2} (\frac{1}{q})_{a_i-b_i}}
\frac{P(b_1,\cdots,b_n)}{P(a_1,\cdots,a_n)}.\] The resulting $n$-tuple
of partitions is distributed according to
$M_{\Gamma,\vec{U}}$. \end{theorem}

	As a final remark, observe that letting $C$ be diagonal with
entries $\frac{1}{P(a_1,\cdots,a_n)}$, one obtains a factorization
$K=CMC^{-1}$, where $M$ is defined by 

\[ M(\vec{a},\vec{b}) = \prod_{1 \leq i \leq j \leq n} q^{f_{ij} a_i
a_j} \prod_{i=1}^n \frac{U_i^{a_i}}{q^{a_i^2}
(\frac{1}{q})_{a_i-b_i}}.\] A very natural problem is to investigate
the eigenvector matrix $E$ of $M$, and $E^{-1}$, in order to obtain
new Bailey Lemmas.

\section{Acknowledgments} This research was supported by an NSF
Postdoctoral Fellowship. The author thanks Persi Diaconis for
discussions and a referee for correcting historical mistakes in the
bibliography.

Current address: Stanford University Math Department, Building 380, MC
2125, Stanford CA 94305, USA.

Current email: fulman@math.stanford.edu

\end{document}